\documentclass[12pt,a4paper]{amsart}
\usepackage{amsfonts,amssymb}
\numberwithin{equation}{section}

 \def\numset#1{{\mathbb #1}}
 \def\setZ{\numset{Z}}

\theoremstyle{plain}
\newtheorem{Th}{Theorem}
\newtheorem*{Lemma}{Lemma}
\newtheorem*{Cor}{Corollary}

 \theoremstyle{definition}
\newtheorem*{Def}{Definition}

\newtheorem*{Sejt}{Conjecture}

\newtheorem{K}[Th]{Problem}

\newcommand{\eq}[1]{\eqref{#1}}

\newcommand\ndiv{\nmid}
\newcommand{\setN}{{\mathbb N}}
\newcommand{\szd}{\triangle}

\begin{document}

\title{Additive decomposition of signed primes}

\author{Imre Z. Ruzsa}
\address{Alfr\'ed R\'enyi Institute of Mathematics\\
     Budapest, Pf. 127\\
     H-1364 Hungary
}
\email{ruzsa@renyi.hu}
 \thanks{Supported by  NKFI grants K-129335, K-119528, KKP-133819.}

     \begin{abstract}
     Under the prime-tuple hypothesis, the set of signed primes is a sumset.
     \end{abstract}

  \maketitle
In memoriam Andreae Schinzel, excellentis mathematici et optimi viri.

     \section{Introduction}
   
Ostmann's ``inverse Goldbach'' problem (this term is probably Wirsing's) asks whether the set
$P$  of (positive) primes is asymptotically a sumset, that is, whether there are sets $A, B$ 
(having more than 1 elements each)
such that  $P \szd( A+B)$ is finite, with $\szd$ meaning  symmetric difference.

Still unsolved, generally a negative answer is expected. See Elsholz\cite{elsholtz01, elsholtz06}
 for the best partial results.

The aim of this paper is to show that the situation changes if we admit negative primes.

  \begin{Th}
    Assuming the { prime-tuple hypothesis} there are infinte sets
$A, B \subset \setZ$ such that $A+B$ is exactly the collection of (positive and negative) primes
satisfying $|p| > 3$; moreover,  every prime $p$ has exactly one representation as 
$p=a+b$, $a\in A$, $b\in B$.
  \end{Th}
  
   It is easily seen that $\pm 2, \pm 3$ must be omitted, and that if the summands have more than 1
elements, they both must be infinite.
  
Some details of the argument are easier written in terms of differences, so we reformulate the
theorem as follows.

\begin{Th}[Variant with differences]
 Assuming the {prime-tuple hypothesis} there are infinite sets
$A, B $ of positive integers such that $A-B$ is exactly the collection of (positive and negative) primes
satisfying $|p| > 3$; moreover,  every prime $p$ has exactly one representaion as 
$p=a-b$, $a\in A$, $b\in B$.
  \end{Th}

\section{The prime-tuple hypothesis}

The prime-tuple hypothesis, generally assumed to be true but hopeless, expresses that
linear forms can simultaneously represent primes unless there is a congruence obstacle.

  \begin{Sejt}[Prime-tuple hypothesis]
    Let $a_i, b_i$ be integers, $a_i \neq 0$. There are infinitely many values of $x$
such that all $a_ix+b_i$ are prime, unless there is a prime $p$ such that for all
 $x$ we have $p | a_ix+b_i$ for some $i$.
  \end{Sejt}

(Simplest case is the twin prime conjectue.)

We shall apply the following special case.

\begin{Cor}
  Given $d_i\in \setZ$, $q>1$ and $t$, the condition for the existence of infinitely many
 $x \equiv t \pmod q$ such that all $x+d_i$ are primes is that for all primes $p$:

--- if $p|q$, then $p \ndiv t+d_i$ for any $i$;

--- if $p\ndiv q$, then the
  $d_i$ do not contain a complete system of residues modulo $p$.
\end{Cor}

This is the preceeding applied to the forms $qx+t+d_i$.

\section{The plan}

 Assume we have finite sets $A,B$ such that all elements of $A-B$ are primes. We want to represent
a further prime $r$. How can we do it? Try

$A'=A \cup \{ x\}$, $ B' = B \cup \{ x-r\}$. 

This works if all elemens of  $B-x$ and $A-x+r$ are prime.

We find such an $x$ if $B \cup (A+r)$ is not a complete system modulo any prime.

But: the inclusion of $x$ and $x-r$ may spoil this property; we may build a trap to kill the plan
after several steps.

Remedy: we will a priori restrict the possible residues
 mod $p$. We clearly cannot do this initially for all primes; we shall dynamically
add more and more restricions, that is, new elements compatible with given restrictions
and new restrictions compatible with given elements.

\section{Prime-compatible sets of residues}

\begin{Def}
  Let $p$ be a prime and $U, V \subset \setZ_p$ nonempty sets of residues. We say that $U,V$
  form a \emph{prime-compatible pair}, if
   \[  (U \setminus V) - (V \setminus U) = \setZ_p^* .\]
\end{Def}

The residues of $A$ and $B$ modulo $p$ must have this property.

\begin{Lemma}
  Let $p$ be a prime.

(a) Let $W\subset\setZ_p$ be a nonempty set of residues. If
\begin{equation}\label{condi}
   |W| < \frac{p-1}{2} -  \frac{\log p}{\log 4/3} ,\end{equation}
then there are prime-compatible sets modulo $p$ satisfying $W\subset U$, $W\subset V$.

(b) For $p \geq7$ there are prime-compatible sets modulo $p$ such that $1,11 \in U\setminus V$, $ 6\in V\setminus U$ and $|U\cap V| =2$.
\end{Lemma}

Part (b) is motivated by the following consideration. By looking at the residues of $A$ and $B$ modulo 5
it is easy to see that the representations of 5 and $-5$ must have a common element; we set $1,11\in A$ and 
$6\in B$. 
\begin{proof}
 (a): Write $|W| =n$. We assign the $p-n$ elements of $\setZ_p \setminus W$ randomly to $U$ or $V$, and estimate the probability
that a $z\in\setZ_p^*$ is represented as a difference of these new elements.

Consider the pairs $(z,2z), (3z,4z), \ldots, ((p-2)z, (p-1)z)$.
From these $(p-1)/2$ pairs at most $n$ contain an element of $W$, at least $ (p-1)/2-n$ remain.

A pair represents $z$ with probability 1/4, so 
 \[ \Pr(\text{none is good }) \leq  (3/4)^{\frac{p-1}{2} - n} .\]

If this quantity is $<1/p$, there is a choice that works for all $z$.
This yields the condition
 \[  (3/4)^{\frac{p-1}{2} - n} < \frac{1}{p}, \]
 which can be rearranged as
  \[ \frac{p-1}{2} - n > \frac{\log p}{\log 4/3} ,\]
that is, \eq{condi}.

(b): Let $U'=\{ 0,1,2, p-1 \}$, $V'=\{ 0,1, 3, 4, \ldots, p-2\}$. This is a prime-compatible pair modulo $p$.
This property is preserved by linear transformations. Define $\alpha,\beta$ by
 \[ 3\alpha \equiv -10 \pmod p, \ \beta=\alpha+11 .\]
The linear function $\alpha x+\beta$ maps 2 to 1, $p-1$ to  11 and $(p+1)/2$ to 6, so the sets $U=\alpha\cdot U'+\beta$, $V=\alpha\cdot V'+\beta$
are suitable.
\end{proof}

The bound could be improved with little effort, but it is irrelevant for our application.

\section{The construction}

We prove the theorems.

 Let  $r_1, r_2, \ldots$ be the sequence of signed primes  $ \pm p$, $|p| > 3$, ordered by increasing absolute value
 ($p_1=5, p_2=-5$, etc.) We shall construct  sequences of sets $A_n$, $B_n$ of nonnegative integers such that
 \[ A_1 \subset A_2 \subset \ldots, \ B_1 \subset B_2 \subset \ldots, \]
 \[ A_n-B_n \supset \{ r_1, r_2, \ldots, r_n\} ,\]
  \[ |A_n| \leq n, \  |B_n| \leq n , \]
   \[   |A_n| > |A_{n-1}|,  \  |B_n| > |B_{n-1}| \text{ infinitely often}, \]
and the numbers $a-b, a\in A_n, b\in B_n$ are all distinct primes. Clearly the sets $A=\bigcup A_n$, $B=\bigcup B_n$ 
will have the properties asserted in Theorem 2.

Let $K$ be a constant such that

\begin{equation}\label{bound}
   2n < \frac{|r_n| -1}{2} -  \frac{\log |r_n|}{\log 4/3} \end{equation}
holds whenever $| r_n | \geq K/2$. Such a constant exists since $|r_n| \sim (n \log n)/2$.

We construct the sets recursively. In Step $n$ we will have sets $A_n$, $B_n$ and sets
 $U_p, V_p \subset \setZ_p$ for all primes $p< \max (K, |r_{n+1}|)$ such that 
\[ A_n \pmod p \subset U_p, \ B_n \pmod p \subset V_p ,\]
and $U_p,V_p$  form a {prime-compatible pair}.

The starting point is $A_2=\{1,11\}$, $B_2=\{6 \}$, $U_2=\{1 \}$, $V_2=\{0 \}$, $U_3=\{1,2 \}$, $V_3=\{0 \}$, 
$U_5= \{0,1,2 \}$, $V_5=\{1,3,4 \}$,
and for $7 \leq p <K$ the sets $U_p, V_p$ are  given in part (b) of the Lemma.

Assume we have $A_{n-1}$, $B_{n-1}$ and $U_p, V_p \subset \setZ_p$  for all primes $p< \max (K,|r_{n}|)$.
We construct $A_n, B_n$ and $U_p, V_p$ for $\max(K, |r_n|) \leq p < \max (K,|r_{n+1}|)$. 
(This is either $p=|r_n|$, or there is no such prime.)

If $r_n\in A_{n-1}-B_{n-1}$, we put  $A_n=A_{n-1}$, $B_n=B_{n-1}$. Otherwise we will set
 \[ A_n=A_{n-1} \cup \{ x\}, \ B_n=B_{n-1}\cup\{x-r_n \} \]
with some positive integer $x$. This $x$ should have the following properties:

(i)  all elements of $x-B_{n-1}$ and $A_{n-1} - x+r_n$ are prime;

(ii) $x \pmod p \in U_p$, $x-r_n  \pmod p  \in V_p$ for  all primes $p< \max (K,|r_{n}|)$;

(iii) no coincidence.

It is clear that condition (iii) excludes only finitely many values of $x$, hence if we
find infinitely many that satisfy (i) and (ii), we are done.

First we fix $x$ mod $p$ for $p< \max(K,|r_n|)$, $p \neq |r_n|$. ($p=|r_n|$ is a possibility for the first
few primes.)

There are  $u_p \in U_p \setminus V_p$, $v_p\in V_p\setminus U_p$, such that $$u_p-v_p \equiv r_n \pmod p . $$
Impose
 \[ x \equiv u_p \pmod p, \ x-r_n \equiv v_p \pmod p.  \]
If $p=|r_n|$, we must have $u_p=v_p$. Now we use the 2 element of $U_p\cap V_p$, choosing different ones for $p$ and $-p$.

According to the Corollary to the prime-tuple hypothesis, the conditions for the existence of 
values of $x$ that satisfy condition (i) above is the following:

--- if $p< \max(K,|r_n|)$, then $p \ndiv u_p-b$ for $b\in B_{n-1}$ and 
$p \ndiv u_p-(a+r_n) = v_p-a$ for $a\in A_{n-1}$;

--- if $p\geq  \max(K,|r_n|)$, then $B_{n-1}\cup (A_{n-1}+r_n)$ does not contain a complete system of residues modulo $p$.

To check the first, note that for $p \neq |r_n|$ we have $b \pmod p \in V_p$ and $u_p \notin V_p$; similarly
$a \pmod p \in U_p$ and $v_p \notin U_p$. For $p=|r_n|$ this means that $A_{n-1}\cup B_{n-1}$ avoids the residue class of $U_p\cap V_p$
assigned to $r_n$. Indeed, in previous stepts for $r_j$ with $|r_j| \neq p$ we asigned elements of $U_p \setminus V_p$ and $V_p\setminus U_p$
to $r_j$, and for $|r_j| = p$ (which may be the case for $j=n-1$) we used the other element of $U_p\cap V_p$.

To check the second, observe that $|B_{n-1}\cup (A_{n-1}+r_n)| <2n $ and  $\max(K,|r_n|) > 2n$ by \eq{bound}.

Finally we construct $U_p, V_p$  modulo $p= |r_n|$ (if necessary). The requirement is that they form a
prime-compatible pair and contain the at most $2n$ residues of $A_n\cup B_n$. The existence of such sets follows from
part (a) of the Lemma.

This ends the proof of the Theorems.

\section{Remarks}

We don't need the full strength of the prime-tuple hypothesis; however, if the primes form a sumset,
then \emph{some} configurations must appear infinitely many times.

\begin{K}
  Can one prove, \emph{without resorting to the prime-tuple hypothesis}, the existence
of infinite sets $A, B \subset\setN$ such that   $A+B\subset P$?
\end{K}

Recent advances related to the twin prime conjecture imly that we can find infinite $A$ and
arbitarily large finite $B$.

Granville\cite{granville90} proved that under the prime-tuple hypothesis the set of primes
contains very general sorts of sumsets.

\begin{K}
  Can one prove, \emph{without resorting to the prime-tuple hypothesis}, the existence
of infinite sets $A, B \subset\setN$ such that no element of  $A+B$ is divisible by any prime of form
 $4k+1$?
\end{K}

For $4k-1$ we have an easy example.

 \bibliographystyle{amsplain}
     \bibliography{cimek}


     \end{document}